\documentclass[12pt,twoside]{article}
\usepackage{tikz}
\usepackage[T2A]{fontenc}
\usepackage[cp1251]{inputenc}
\usepackage[english]{babel}
\usepackage{floatflt}
\usepackage{listings}

\usepackage{amssymb}
\usepackage{amsfonts}
\usepackage{amsmath}
\usepackage{amsthm}
\usepackage{graphicx,graphics}
\usepackage{caption,psfrag}
\usepackage{xspace}
\usepackage{comment}

\theoremstyle{plain}
\newtheorem{lemma}{Lemma}
\newtheorem{theorem}{Theorem}
\newtheorem{definition}{Definition}

\newtheorem{proposition}{Proposition}
\newtheorem{remark}{Remark}

\newcommand{\myRe}{\operatorname{Re}}
\newcommand{\myIm}{\operatorname{Im}}
\newcommand{\sign}{\operatorname{sign}}
\newcommand{\FlatOneSign}{\operatorname{sign^f_1}}
\newcommand{\FlatTwoSign}{\operatorname{sign^f_2}}

\newcommand{\ring}{\mathcal{R}}
\newcommand{\algA}{\mathcal{A}}

\newcommand{\barbar}[1]{\bar{\bar{#1}}}

\newcommand{\Input}{{\operatorname{inp}}}
\newcommand{\Output}{{\operatorname{out}}}

\newcommand{\Left}[1]{{\operatorname{L}(#1)}}
\newcommand{\Right}[1]{{\operatorname{R}(#1)}}

\begin{document}

\author{Vladimir Tarkaev
\footnote{
The work is supported by RSF (grant number 23-21-10014).
}
}

\title{
An analogue of Turaev comultiplication for knots
in non-orientable thickening of a non-orientable surface
}

\maketitle

\thispagestyle{empty}
\vspace{3mm}
\noindent{\small{\it Chelyabinsk State University, Chelyabinsk, Russia}}

\noindent{\small{\it Krasovskii Institute of Mathematics and Mechanics, Ural Branch of the Russian Academy of Sciences, Yekaterinburg, Russia
}}

\noindent{\small{\it v.tarkaev@mail.ru}}

\begin{abstract}
This paper concerns pseudo-classical knots in the non-orientable manifold $\hat{\Sigma} =\Sigma \times [0,1]$,
where $\Sigma$ is a non-orientable surface
and a knot $K \subset \hat{\Sigma}$ is called pseudo-classical if $K$ is orientation-preserving path in $\hat{\Sigma}$.
For this kind of knot we introduce an invariant $\Delta$
that is an analogue of Turaev comultiplication for knots in a thickened orientable surface.
As its classical prototype, $\Delta$ takes value in a polynomial algebra generated by homotopy classes of non-contractible loops on $\Sigma$,
however, as a ground ring we use some subring of $\mathbb{C}$ instead of $\mathbb{Z}$.
Then we define a few homotopy, homology and polynomial invariants,
which are consequences of $\Delta$,
including an analogue of the affine index polynomial.

Keywords:
knots in non-orientable manifold,knots in thickened surface,invariants of knots,Turaev comultiplication,affine index polynomial.
\end{abstract}

%
\markboth{}{} \markboth{Vladimir Tarkaev
}{\footnotesize An analogue of Turaev comultiplication}

\section*{Introduction}

Let $\Sigma$ be a compact, non-orientable surface
that is not necessarily closed.
In this work, we deal with knots in the non-orientable manifold
$\hat{\Sigma} =\Sigma \times [0,1]$.
Under an oriented knot $K \subset \hat{\Sigma}$, we mean the image of a smooth  embedding of an oriented circle into the interior of $\hat{\Sigma}$.
Knots are considered up to ambient isotopy of $\hat{\Sigma}$.
Likewise, in the case of knots in a thickened orientable surface, 
knots in $\hat{\Sigma}$ can be represented by diagrams
obtained via an appropriate  projection map $p: \hat{\Sigma} \to \Sigma \times \{ 0 \}$.
As usual, a diagram of a knot is a $4$-valent graph embedded into the surface
whose vertices are equipped with standard over/under information.
Two knot diagrams on $\Sigma$ represent the same knot in $\hat{\Sigma}$ if and only if 
one can be transformed into the other by a finite sequence consisting of isotopies of $\Sigma$
and Reidemeister moves, which are completely analogous to classical ones.
As in the orientable case,
two generic loops on $\Sigma$ are homotopic
if and only if one can be transformed into the other by a finite sequence consisting of
ambient isotopies and flat versions of Reidemeister moves.

Now, we briefly review the construction of the affine index polynomial
of knots in a thickened orientable surface
\cite{Kauffman2013},\cite{KauffmanLinkingNumber2013}.
Later, we will modify it to work in the non-orientable case.

Let $D$ be an oriented  diagram of an oriented knot in a thickened orientable surface.
Consider a crossing $x \in \#D$, where $\#D$ denotes the set of crossings of the diagram $D$.
We associate with $x$ two closed directed paths $l_1(x)$ and $l_2(x)$ in $D$.
Both of these paths start and end at $x$ and travel along $D$ in the positive direction:
$l_1(x)$ goes away from $x$ by the upper outgoing arc and returns to $x$ by the lower incoming one,
while $l_2(x)$, on the contrary,  goes away from $x$ by the lower outgoing arc and returns by the incoming upper one.
The affine index polynomial can be defined by the following equality:
$$
P(K) =\sum_{x \in \#D}\sign(x) \big(t^{l_1(x) \cdot l_2(x)} -1 \big) \in \mathbb{Z}[t,t^{-1}],
\eqno(1)
$$
where $D$ is a diagram of the knot $K$,
$\sign(x)$ denotes the sign of the crossing $x$
and $l_1(x) \cdot l_2(x)$ denotes the intersection number of $l_1(x),l_2(x)$ viewed as directed closed curves on the corresponding surface.

This polynomial invariant can be viewed as a consequence of
Turaev comultiplication $\Delta$
\cite{Turaev1991},\cite{Turaev2004}
that in our terms can be written in the form
$$
\Delta(K) =\sum_{x \in \#D}\sign(x) [l_1(x)] \otimes [l_2(x)] \in \mathcal{A},
\eqno(2)
$$
where$D$ is a diagram of an oriented knot $K$,
$[l_j(x)],j=1,2,$ denotes the free homotopy class of specified loop and
$\mathcal{A}$ is the polynomial $\mathbb{Z}$-algebra generated by
free homotopy classes of non-contractible loops on the underlying surface.

The definitions~(1) and~(2) 
do not work for diagrams of knots in the manifold $\hat{\Sigma}$
because they use  two classical notions (the sign of a crossing and the intersection number of loops),
which are defined via the orientation of the corresponding surface.
Note that knots in the orientable $I$-bundle over a non-orientable surface
can be represented by diagrams in appropriate polygon
whose sides are assumed to be identified (maybe) with a twist
~\cite{Bourgoin2008}.
In this case, the classical concept of the sign of a crossing works as well,
however, it is so only because of orientability of the corresponding manifold.

In our case, when both of the surface and its thickening are non-orientable,
the classical definition of the sign of a crossing is not applicable.

In this paper we propose a definition of the sign of a crossing
that does not need the concept of the orientation,
and in the classical situation it is equivalent to the classical one.
In the non-orientable case, the approach works for diagrams
representing knots that viewed as loops in $\hat{\Sigma}$
are orientation-preserving paths.
We call them \emph{pseudo-classical knots}
(see Section~1).
Our key idea is to consider $2$-cabling of the diagram in question
and to choose a value of the sign of a crossing by looking at the $4$-crossings pattern corresponding to the crossing of the initial diagram.
Precise definitions are given in Section~2.
In Section~3
using the notion of the sign of a crossing,
we define analogues of the linking number for pseudo-classical links
 and analogues of the intersection number and Goldman bracket
for pseudo-classical curves.

In Section~4
we define an analogue of Turaev comultiplication for pseudo-classical knots,
and then in Section~5
we define some invariants,
which are consequences of this analogue of Turaev comultiplication,
including an analogue of the affine index polynomial.

\section{Knots and diagrams}
\label{sec:Knots}

Throughout  $\Sigma$ denotes a non-orientable surface with (maybe) non-empty boundary
that is not necessarily closed,
and $\hat{\Sigma}$ denotes the non-orientable manifold of the form $\Sigma \times [0,1]$
with fixed structure of the direct product.
The latter condition becomes crucial for us in the case
of surface with non-empty boundary  
when there are aforementioned structures of the corresponding manifold
with homeomorphic but non-isotopic base surfaces.

Knots/links in $\hat{\Sigma}$, generic curves and knot/link diagrams in $\Sigma$ can be defined analogously to the case
of a thickened orientable surface
(below we will refer to the latter case  as the ``classical case'').
Knots in $\hat{\Sigma}$ will be regarded up to ambient isotopy.

One can show that in our case, likewise, in the classical one,
two diagrams on $\Sigma$ represent the same knot
if and only if one can be transformed into the other  by a finite sequence of Reidemeister moves $R_1,R_2,R_3$
(which are just the same as classical ones)
and ambient isotopies.

A knot $K \subset \hat{\Sigma}$ is called a \emph{pseudo-classical}
if it is an orientation-preserving path in $\hat{\Sigma}$.
In other words, a knot in $\hat{\Sigma}$ is pseudo-classical
if its regular neighbourhood in the manifold is homeomorphic to the solid torus.
Clearly, in non-orientable manifolds not all knots have the property,
since in this case there are knots whose regular neighbourhood is homeomorphic to the solid Klein bottle.

Let $K$ be a pseudo-classical oriented knot represented by a diagram $D \subset \Sigma$.
Since the knot is pseudo-classical,
the diagram (or, more precisely, the projection of $K$ onto $\Sigma$)
viewed as a closed directed path on the surface
is an orientation-preserving  path.
(However, it can contain loops
that are not orientation-preserving paths).
Hence $2$-cabling $D^2 \subset \Sigma$ of $D$
is a diagram of a $2$-component link
(recall that $2$-cabling of a diagram is, informally speaking, the same diagram
drawn by doubled-line with preserving over/under information in all the crossings).
Note, again, that this is not the case for an arbitrary knot diagram in $\Sigma$.
That is because $\Sigma$ is assumed to be non-orientable, hence there are knot diagrams in the surface
which are orientation-reversing paths,
hence their $2$-cabling represents not a $2$-component link but a knot
which goes along the given knot twice.
In the latter case, $2$-cabling bounds the {M}\"obius strip (going across itself near each crossing of the given diagram)
while $2$-cabling of a pseudo-classical knot bounds an annulus (having self-intersections).

In the case of an oriented surface dealing with $2$-cabling of an oriented knot diagram, one can speak of left and right components.
Since $\Sigma$ is non-orientable, we do not have the concepts of left and right,
 however, below we will use the terms ``left component'' and ``right component'' of $D^2$
keeping in mind that in our situation this is nothing more than a convenient notation.
A choice (which component is chosen as the left and which is chosen as the right) will be called the \emph{labeling} of $D^2$ (or labeling of $D$ for simplicity).
The components of $D^2$ will be denoted by $\Left{D}$ and $\Right{D}$, respectively.
Throughout, $\Left{D}$ and $\Right{D}$ are assumed to be oriented
so that their orientations are agreed with the orientation of $D$.
Sometimes we will rename $\Right{D}$ into $\Left{D}$ and vice versa,
the transformation will be called the \emph{relabeling} of $D$.

A single crossing of $D$ becomes in $D^2$ a pattern of $4$ crossings
(see Fig.~1).
There are two distinguished crossings  in each of these patterns:
the first one in which components of $D^2$
come into the pattern and the second one in which they go out.
These crossings are called the \emph{input} (resp. \emph{output})  crossing for the crossing $x$ 
and they are denoted by $\Input(x)$ (resp. $\Output(x)$\,),
where  $x \in \#D$ is that crossing of $D$ to which the pattern corresponds.
Here and below, we denote by $\#D$ the set of all the crossings of a diagram $D$.

Note if the surface under consideration is orientable
then input and output crossings are necessarily intersections of distinct components of $2$-cabling
while in non-orientable case
the situation when input and output crossings are a self-crossing of the same component of $D^2$ may occur.

\begin{figure}
\centering 
$$
\begin{array}{cc}
\includegraphics[width=0.45\textwidth]{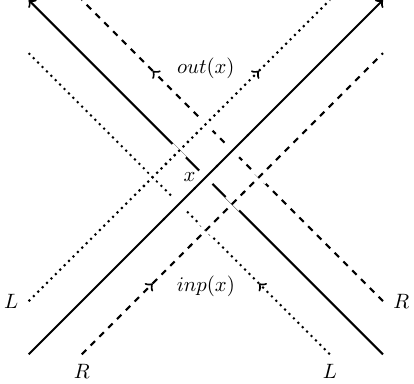} &
\includegraphics[width=0.45\textwidth]{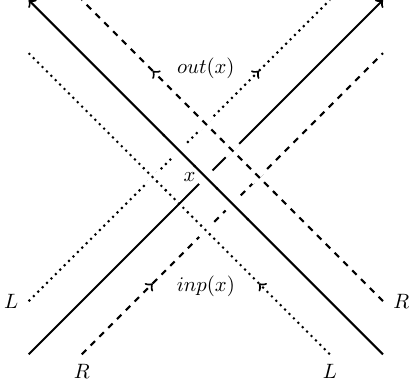} \\
(a) & (b) \\
\includegraphics[width=0.45\textwidth]{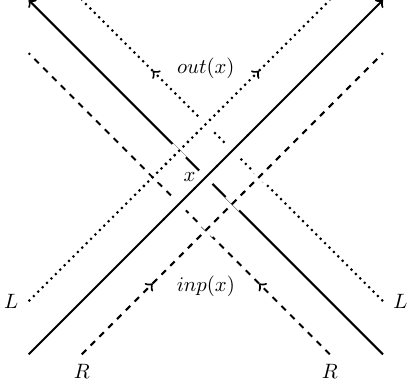} &
\includegraphics[width=0.45\textwidth]{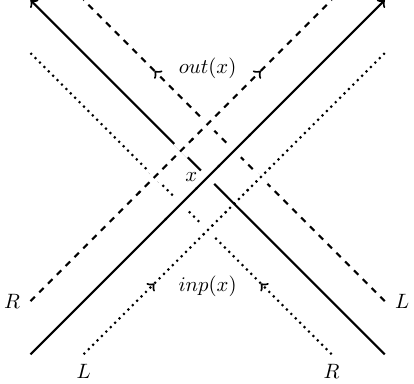} \\
(c) & (d)
\end{array}
$$
Fig. 1: Segments belong to $D, \Right{D}$ and $\Left{D}$ are drawn by solid, dashed and dotted lines, respectively.
\label{fig:1}
\end{figure}

\section{A definition of the sign of a crossing}
\label{sec:Sign}

\subsection{The sign of a crossing equipped with over/under information}
\label{sec:SignOfDiagramCrossing}

Recall that in the classical situation the sign of a crossing $x$ of an oriented diagram $D$ on an oriented surface $F$
is defined to be $1$ (resp. $-1$) if the pair $(t_1, t_2)$
is positive (resp. negative) basis
in the tangent space  of $F$ at $x$,
where $t_1,t_2$ are positive tangent vectors of overgoing and undergoing branches of $D$ at $x$.
Clearly, the  definition has no sense if the surface  is non-orientable.
In the following definition, we use a labeled $2$-cabling of the diagram  under consideration instead of the orientation of the surface,
in which the diagram lies.

\begin{definition}
\label{def:Sign}
Let $D \subset \Sigma$ be a diagram of an oriented pseudo-classical knot and $D^2$ be a $2$-cabling of $D$.
The \emph{sign} of a crossing $x \in \#D$ is defined to be the mapping
$\sign : \#D \to \mathbb{C}$
given by the following rule:
\begin{equation*}
\label{eq:Sign}
\sign(x) =\begin{cases}
1 & \text{if $\Right{D}$ goes over $\Left{D}$ in $\Input(x)$ \textnormal{(Fig.~1(a))}),} \\
-1 & \text{if $\Left{D}$ goes over $\Right{D}$ in $\Input(x)$ \textnormal{(Fig.~1(b))}),} \\
i & \text{if $\Right{D}$ goes over itself in $\Input(x)$ \textnormal{(Fig.~1(c))}),} \\
-i & \text{if $\Left{D}$ goes over itself in $\Input(x)$ \textnormal{(Fig.~1(d))}),}
\end{cases}
\end{equation*}
where $i \in \mathbb{C}$ is the imaginary unit. 
\end{definition}

In the classical situation, when both the surface and its thickening are orientable,
it is natural to choose as $\Right{D}$ that component of $D^2$ which goes to the right of $D$ with respect to its actual orientation,
then the above rule  gives the standard sign of a crossing.

Note that a crossing $x \in \#D$ decomposes $D$ into two closed paths $l_1(x),l_2(x)$.
Both of these paths start and end at $x$ and travel along $D$ in the positive direction:
$l_1(x)$ goes away from $x$ by the upper outgoing arc and returns to $x$ by the lower incoming one,
while $l_2(x)$, on the contrary,  goes away from $x$ by the lower outgoing arc and returns by the incoming upper one.
A pseudo-classical  diagram $D$ viewed as a directed path on $\Sigma$
is an orientation-preserving path.
Hence, if $\sign(x) =\pm 1$ 
then both of $l_1(x)$ and $l_2(x)$ are paths preserving the orientation,
while if $\sign(x) =\pm i$ then each of $l_1(x)$ and $l_2(x)$ is a path  reversing the orientation.
The situation when one of $l_1(x),l_2(x)$ preserves the orientation while the other does not
is impossible in our context.

The following proposition is a direct consequence  of our definition of the sign
and the observation that the input crossing of each $4$-crossing pattern becomes the output crossing and vice versa
as a result of the reversing orientation of $D$.
Recall that the switching operation in a crossing of a diagram
consists in permuting roles of the overgoing and the undergoing arcs in the crossing,
i.e., in reassigning the overgoing arc into the undergoing one and vice versa.

\begin{proposition}
\label{proposition:Sign}
Let $D \subset \Sigma$ be a diagram of an oriented pseudo-classical knot
and $x \in \#D$.
Then:

{\bf 1}.
The relabeling of $D$ transforms $\sign(x)$ into $-\sign(x)$.

{\bf 2}.
The reversing of the orientation of $K$ (and consequently $D$)
while the labeling of $D$ remains unchanged
transforms $\sign(x)$ into $-\sign(x)$.

{\bf 3}.
The switching in $x$ transforms $\sign(x)$ into
$-\myRe \sign(x) +i \myIm \sign(x) =-\overline{\sign(x)}$.
\end{proposition}

\subsection{The sign of a crossing in a diagram of pseudo-classical link}
\label{sec:SignForLinks}

Using the same idea, we define the sign of a crossing
in a diagram of a pseudo-classical link,
i.e., a link each of which component is a pseudo-classical path in $\hat{\Sigma}$.
Namely, consider a link $\ell =\ell_1 \amalg \ldots \amalg \ell_m \subset \hat{\Sigma}$
so that $\ell_j$ is an orientation-preserving path in $\hat{\Sigma}$ for any $j=1,\ldots,m$,
and let $D, D^2 \subset \Sigma$ be a diagram of $\ell$ and its $2$-cabling, respectively.
The components of $D$ and their $2$-cablings will be denoted by $D_j$ and $D_j^2$, $j=1,\ldots,m$, respectively.

\begin{definition}
\label{def:SignForLink}
Let $D$ and $D^2$ be as above.
If $x \in \#D$ is an intersection point of $D_i$ and $D_j$, $i \not=j$,
Then the sign of the crossing is defined to be
$$
\sign(x) =\begin{cases}
1 & \text{if $\Right{D_i}$ (resp. $\Right{D_j}$) goes over $\Left{D_j}$ (resp. $\Left{D_i}$) in $\Input(x)$,} \\
-1 & \text{if $\Left{D_i}$ (resp. $\Left{D_j}$) goes over $\Right{D_j}$ (resp. $\Right{D_i}$) in $\Input(x)$,} \\
i & \text{if $\Right{D_i}$ (resp. $\Right{D_j}$) goes over $\Right{D_j}$ (resp. $\Right{D_i}$) in $\Input(x)$,} \\
-i & \text{if $\Left{D_i}$ (resp. $\Left{D_j}$) goes over $\Left{D_j}$ (resp. $\Left{D_i}$) in $\Input(x)$,}
\end{cases}
$$
where $i \in \mathbb{C}$ is the imaginary unit. 

If $x \in \#D$ is a self-intersection of a component of $D$,
then $\sign(x)$ is understood in the sense of Definition~\ref{def:Sign}.
\end{definition}

If $x \in \#D$ is a self-intersection point of a component,
then the effect of the relabeling, reversing and switching
is the same as in the case of knot diagram.
While, if a crossing is an intersection point of two distinct components,
the labeling and the direction of both of the intersecting components affect the resulting value of the sign.
Hence, a proposition analogous to Proposition~\ref{proposition:Sign}
for such a crossing takes the form:

\begin{proposition}
\label{proposition:SignForLink}
Let $D \subset \Sigma$ be a diagram of a pseudo-classical link.
If $x \in \#D$ is an intersection point of components $D_i$ and $D_j$, $i \not=j$,
in which $D_i$ goes over $D_j$,
Then The following transformations of $D$ transform $\sign(x)$ into:
$$
\begin{array}{rcl}
\textnormal{Transformation} & & \textnormal{The resulting value} \\
\text{relabeling of } D_i & \quad & S(-\sign(x)), \\
\text{relabeling of } D_j & \quad & S(\sign(x)), \\
\text{simultaneous relabeling of } D_i, D_j & \quad & -\sign(x), \\
\text{reversing of } D_i & \quad & S( \sign(x)), \\
\text{reversing of } D_j & \quad & S(-\sign(x)), \\
\text{simultaneous  reversing of } D_i,D_j & \quad & -\sign(x), \\
\text{switching in } x & \quad & -\myRe \sign(x) +i \myIm \sign(x),
\end{array}
$$
where $S: \mathbb{C} \to \mathbb{C}$ denotes
the symmetry along the bisectrix of the first and the third coordinate angles.
\end{proposition}
Note that $S(-z)$ is
the symmetry along the bisectrix of the second and the fourth coordinate angles.

\begin{remark}
\label{remark:Bigon}
Let crossings $x_1$ and $x_2$ be vertices of a bigon,
which appears as a result of the second Reidemeister move,
i.e., the same branch goes over in both of these crossings.
In the classical case, in this situation $\sign(x_1) =-\sign(x_2)$.
It is easy to check that the same equality holds for diagrams of pseudo-classical knots and pseudo-classical links
independently on whether involved branches of the diagram have the same direction (Fig.~2 on the left)
or not (Fig.~2 on the right).

\begin{figure}[h!]
\centering 
$$
\begin{array}{cc}
\includegraphics[width=0.45\textwidth]{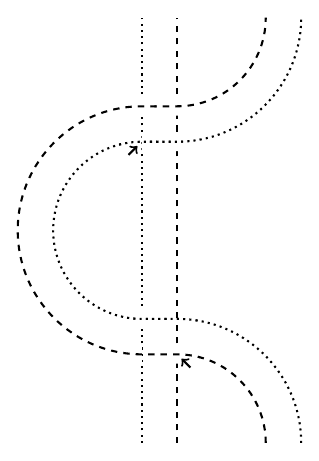} &
\includegraphics[width=0.45\textwidth]{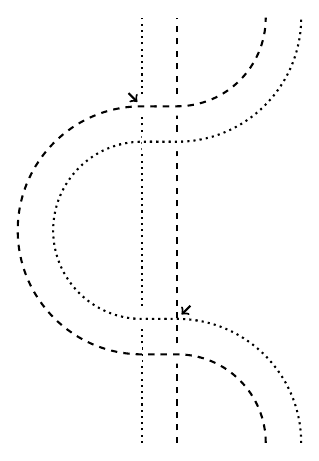}
\end{array}
$$
Fig. 2:
\small{
$\Right{D}$ and $\Left{D}$ are drawn by dashed and dotted lines, respectively.
Bent strand is directed from bottom to top on both pictures,
while vertical strand is directed from bottom to top on the left and from top to bottom on the right.
Short arrows point to input crossings of $4$-crossings patterns.
}
\end{figure}
\end{remark}

\subsection{The sign of a flat intersection point}
\label{sec:SignOfFlatIntersection}

Under a \emph{flat} intersection point, we understand
a transverse intersection point of two generic curves on $\Sigma$ (or a self-intersection of a generic curve)
that is not equipped with over/under information.
The sign of a flat intersection point can be defined directly via the actual labeling of the intersecting curves
in the same manner as in Definition~\ref{def:Sign}.
But the case of our main interest involves 
flat intersection points obtained
by forgetting over/under information
in crossings of a pseudo-classical diagram.
That is because we will define the sign of a flat intersection point via the sign of non-flat crossing.

\subsubsection{The sign of a flat intersection point of distinct pseudo-classical curves}
\label{sec:SignOfFlatForCurves}

Let $C_1,C_2 \subset \Sigma$ be closed oriented generic pseudo-classical curves.
Let, as usual, $C_1 \cap C_2$ be finite
and each $x \in C_1 \cap C_2$ is a transversal intersection of exactly two branches.
Now, we define $\FlatTwoSign(x)$ --- the sign of a flat intersection point $x \in C_1 \cap C_2$
of an ordered pair of curves
using the sign of non-flat crossing that we defined in the previous section.
The superscript ``f'' in $\FlatTwoSign(x)$ stands for ``flat''
and the subscript ``$2$'' means that $x$ is the intersection point of two distinct curves.

We define $\FlatTwoSign(x)$ to be equal $\sign(x)$ 
in the case when $x$ is equipped with over/under information,
so that $C_1$ goes over $C_2$ in the crossing.
Clearly, $\FlatTwoSign(x)$ depends both on the actual labeling of the curves and on which curve is the first.

Directly from the above definition and Definition~\ref{def:Sign}
we have the following properties:

\begin{proposition}
\label{proposition:FlatSign}
Let $C_1, C_2 \subset \Sigma$ be curves as above
and $x \in C_1 \cap C_2$.
Then:

{\bf 1}.
The following transformations transform $\FlatTwoSign(x)$ into:
$$
\begin{array}{rcl}
\textnormal{Transformation} & & \textnormal{The resulting value} \\
\text{relabeling of $C_1$} & \quad & S(-\FlatTwoSign(x)), \\
\text{relabeling of $C_2$} & \quad & S(\FlatTwoSign(x)), \\ 
\text{simultaneous relabeling of $C_1,C_2$} & \quad & -\FlatTwoSign(x), \\
\text{permuting $C_1$ and $C_2$} & \quad & -\overline{\FlatTwoSign(x)},
\end{array}
$$
where $S: \mathbb{C} \to \mathbb{C}$ denotes
the symmetry along the bisectrix of the first and the third coordinate angles
and $\overline{\FlatTwoSign(x)}$
denotes the complex conjugation.

{\bf 2}.
If we forget over/under information in a crossing $x$ of a pseudo-classical diagram
then 
$$
\FlatTwoSign(x) =\begin{cases}
-\sign(x) & \text{if $C_1$ goes under $C_2$ in $x$ and $\sign(x) =\pm 1$,} \\
\sign(x) & \text{otherwise.}
\end{cases}
$$
\end{proposition}

\begin{remark}
\label{remark:FlatBigonForCurves}
Let $x_1,x_2$ be vertices of  a flat bigon
and let $x_1,x_2$ be intersection points of distinct curves.
Then, by Remark~\ref{remark:Bigon}
and Proposition~\ref{proposition:FlatSign}(2),
$$
\FlatTwoSign(x_1) =-\FlatTwoSign(x_2).
$$
\end{remark}

\subsubsection{The sign of a flat self-intersection point of a pseudo-classical curve}
\label{sec:SignOfFlatForSelfIntersection}

Let $D$ be a pseudo-classical knot diagram.
Consider $x \in \#D$ as a flat self-intersection of the curve
obtained from $D$ as a result of forgetting over/under information in all crossings.

Set
\begin{equation*}
\label{eq:FlatSignSelfIntersection}
\FlatOneSign(x) =\begin{cases}
\sign(x) & \text{if $\sign(x) =\pm i$,} \\
1 & \text{if $sign(x) =\pm 1$.}
\end{cases}
\end{equation*}
Here the superscript ``f'' stands for ``flat''
and the subscript ``$1$'' means that $x$ is a point of self-intersection of a curve.

Directly from the above definition and Proposition~\ref{proposition:Sign}
we have the following:

\begin{proposition}
\label{proposition:FlatSignSelfIntersection}
Let $x$ be a flat self-intersection point of a pseudo-classical curve.
Then, relabeling transforms $\FlatOneSign(x)$ into $\overline{\FlatOneSign(x)}$.
\end{proposition}

\begin{remark}
\label{remark:FlatBigonSelfIntersection}
Let $x_1,x_2$ be the vertices of a bigon
and let both of these points be flat self-intersection of a pseudo-classical curve,
then
$$
\FlatOneSign(x_1) =\overline{\FlatOneSign(x_2)}.
$$
\end{remark}

\section{Invariants of pairs of pseudo-classical curves}
\label{sec:Pairs}

\subsection{An analogue of the linking number}
\label{sec:LinkingNumber}

Consider an oriented diagram $D \subset \Sigma$ of a $2$-component pseudo-classical link,
i.e., $D$ consists of two circles immersed into $\Sigma$
with all standard conditions satisfied.
Let $C_1,C_2$ denote the components of $D$.

Set
$$
l(C_1,C_2) =
\sum_{x \in C_1 \cap C_2}\varepsilon(x) \sign(x),
$$
where
$$
\varepsilon(x)=\begin{cases}
1 & \text{if $C_1$ goes over $C_2$ in $x$,} \\
-1 & \text{otherwise.}
\end{cases}
$$

\begin{proposition}
\label{proposition:LinkingNumber}
The value $l(C_1,C_2)$
is an invariant under Reidemeister moves
of unordered pair $(C_1,C_2)$
up to the following transformations:
\begin{itemize}
\item the symmetry along the bisectrix of the first and the third coordinate angles,
\item the symmetry along the bisectrix of the second and the fourth coordinate angles,
\item multiplying by $-1$.
\end{itemize}
\end{proposition}

\emph{Proof}.
Let the labelings of $C_1$ and $C_2$ be fixed.
Then the invariance of $l(C_1,C_2)$ under Reidemeister moves is a consequence of the following observations:
\begin{itemize}
\item[$R_1$:] Self-intersections of components are not involved in the sum.
\item[$R_2$:] The signs of two involved crossings are opposite (see Remark~\ref{remark:Bigon}).
\item[$R_3$:] The signs of all three crossings involved in the move are preserved.
\end{itemize}
Since there are no canonical labelings of the curves,
we need to consider all four possible labelings.
 The corresponding values can be transformed one into another
by transformations listed in Proposition~\ref{proposition:SignForLink},
which coincide with the ones listed in proving proposition.

If the labelings of the curves are fixed,
then, by definition, permuting of $C_1$ and $C_2$ implies to multiplying $\varepsilon(x)$ by $-1$ for all the intersection points.
At the same time, the value of $l(C_2,C_1)$
coincides with $l(C_1,C_2)$ in the case
of simultaneous relabeling of $C_1$ and $C_2$ (Proposition~\ref{proposition:SignForLink}),
hence the order of the curves in the pair does not matter.
\qed

\subsection{An analogue of the intersection number of pseudo-classical curves}
\label{sec:IntersectionNumber}

Let $C_1,C_2 \subset \Sigma$ be as in previous section. 
Since the curves are oriented and pseudo-classical,
we can use our notion of the sign for their crossings
(see Section~\ref{sec:SignOfFlatForCurves})
to define an analogue of the intersection number.

Set
$$
C_1 \cdot C_2 =\sum_{x \in C_1 \cap C_2}\FlatTwoSign(x).
\eqno(3)
$$

Directly from the above definition and Proposition~\ref{proposition:FlatSign}(1)
we have the following properties:

\begin{proposition}
\label{proposition:IntersectionNumber}
The value $C_1 \cdot C_2$
is an invariant of the ordered pair $(C_1,C_2)$
up to the transformations listed in Proposition~\ref{proposition:FlatSign}(1).
\end{proposition}
In particular, it follows
that this analogue of the intersection number indeed depend on the  order of the intersecting curves
if and only if $\myRe C_1 \cdot C_2 \not=0$ and $\myIm C_1 \cdot C_2 \not=0$.

\begin{remark}
\label{rem:IntersectionNumberForSubloops}
Below, we will use concept of the intersection number
in the following particular  case.
Let $D$ be a diagram of oriented pseudo-classical knot
and let $l_1,l_2 \subset D$ be loops in $D$
whose orientation agrees with the orientation of $D$.
In this case, each intersection point $x \in l_1 \cap l_2$
is a crossing of $D$.
Hence, it can be equipped with $\sign(x)$ coming from a labeling of $D$.
Then we can transform $\sign(x)$ into $\FlatTwoSign(x)$
(see Proposition~\ref{proposition:FlatSign}(2))
and define $l_1 \cdot l_2$ by~(3).
Note that the labelings of $l_1,l_2$ are not independent
because  both of them are inherited from the same labeling of $D$.
Hence, $l_1 \cdot l_2$ is uniquely determined by the actual labeling of $D$
and, by Proposition~\ref{proposition:FlatSign}(1),
relabeling of $D$ implies multiplying $l_1 \cdot l_2$ by $-1$.
Another important difference from the general case
is that $l_1 \cdot l_2$ can be defined independently on whether
the loops $l_1$ and $l_2$ are pseudo-classical or not,
because $\FlatTwoSign$ of their intersection points is determined by $\sign$ of the crossing coming from the diagram.
\end{remark}

\subsection{An analogue of Goldman bracket}
\label{sec:GoldmanBracket}

Let $\ring$ denote the ring
whose elements are complex numbers of the form $a +ib$,
where $a,b \in \mathbb{Z}$ and $i$ is the imaginary unit.
We denote by $\algA$ the polynomial algebra over the ring $\ring$
generated by homotopy classes of oriented closed non-contractible  curve on $\Sigma$
(here we allow curves that are not orientation-preserving paths as well).
Let $\algA_n$ denote the homogeneous degree $n$ component of $\algA$.

Let $C_1,C_2$ be as in previous sections.

We define an analogue of Goldman bracket~\cite{Goldman1986} as follows
\begin{equation*}
\label{eq:GoldmanBracket}
[C_1,C_2] =\sum_{x \in C_1 \cap C_2}\FlatTwoSign(x) [C(x)] \in \algA_1,
\end{equation*}
where $[\cdot]$ denotes the homotopy class of specified loop and
$C(x) \subset \Sigma$ denotes the oriented loop on $\Sigma$
obtained as the result of the smoothing of $x$
that gives the loop whose orientation everywhere agree with the ones of $C_1$ and $C_2$;
some authors call such a smoothing as ``smoothing along the orientation''.

The invariance $[C_1,C_2]$ under flat versions of Reidemeister moves
follows from:

1. The move $R_1$ creates/removes a point of self-intersection
hence it does not affect the value.

2. If the move $R_2$ creates/removes points $x_1,x_2$ so that
$x_1,x_2 \in C_1 \cap C_2$
then (see Remark~\ref{remark:FlatBigonForCurves})
$\FlatTwoSign(x_1) =-\FlatTwoSign(x_2)$
and $[C(x_1)] =[C(x_2)]$.
Hence, the terms corresponding to $x_1$ and $x_2$ cancel out.
If the move $R_2$ creates/removes two points of self-intersection of $C_1$ or $C_2$,
then it does not affect the resulting value.

3. The move $R_3$ preserves terms corresponding to all three involved intersection points.

The changes of $[C_1,C_2]$ as a result of relabeling
are listed in proposition~\ref{proposition:FlatSign}(1).
Therefore, we see that our analogue of Goldman bracket is not a Lie bracket, unlike its classical prototype.

\begin{remark}
\label{rem:GoldmanBracketForSubloops}
As in Remark~\ref{rem:IntersectionNumberForSubloops},
one can consider $[l_1,l_2]$ for any (not necessarily pseudo-classical)
loops in a pseudo-classical diagram.
Relabeling of $D$ transforms $[l_1,l_2]$ into $-[l_1,l_2]$
while $[l_2,l_1] =-\overline{[l_1,l_2]}$.
\end{remark}

\section{An analogue of Turaev comultiplication}
\label{sec:Delta}

\subsection{An analogue of Turaev comultiplication for pseudo-classical knots}
\label{sec:DeltaForKnots}

Let $D \subset \Sigma$ be an oriented diagram of an oriented pseudo-classical knot $K \subset \hat{\Sigma}$.
Now, we can define $\Delta(K)$ by~(2)
in which $\sign(x)$ is understood in the sense
of Section~\ref{sec:SignOfDiagramCrossing}.

\begin{theorem}
\label{theorem:Delta}
Let $D \subset \Sigma$ be a diagram of an oriented pseudo-classical knot $K \subset \hat{\Sigma}$,
then
$$
\Delta{K} =\sum_{x \in \#D}\sign(x) [l_1(x)] \otimes [l_2(x)] \in \algA_2
$$
is an isotopy invariant of the knot $K$
up to multiplying by $-1$.
\end{theorem}

\emph{Proof}.

Fix a labeling of the diagram $D$.

It is sufficient to check
that $\Delta$ is preserved under Reidemeister moves.
Note that Reidemeister moves do not change the homotopy class of loops
$l_j(x),j=1,2,$ if the crossing $x$ lies outside the disk
inside which the move is performed.
Hence, the terms on the right-hand side of~(2)
corresponding to these crossings remain unchanged.

If $x$ is the crossing involved in the move $R_1$,
then either $l_1(x)$ or $l_2(x)$ is contractible
hence $[l_1(x)] \otimes [l_2(x)] =0 \in \algA$.

The move $R_2$ creates/removes a bigon.
Let $x_1,x_2$ be the crossings involved in the move.
Then $[l_j(x_1)] =[l_j(x_2)],j=1,2,$
and (see Remark~\ref{remark:Bigon})
$\sign(x_1) =-\sign(x_2)$.
Hence 
$$
\sign(x_1) [l_1(x_1)] \otimes [l_2(x_1)] +\sign(x_2) [l_1(x_2)] \otimes [l_2(x_2)] =0 \in \algA.
$$

Let $x_1,x_2,x_3$ be crossings involved in the move $R_3$.
Then for any $k=1,2,3,$ the move preserves 
$\sign(x_k)$,
and the homotopy class of the loops
$l_j(x_k),j=1,2,$.
Therefore, $\Delta(D)$ is preserved under the move $R_3$.

Since, in the non-orientable case,
there is no canonical labeling of a diagram,
we need to consider  a relabeling of $D$
which (see Proposition~\ref{proposition:Sign}(1)\,)
multiplies the signs of all the crossings by $-1$,
hence $\Delta(K)$ is an invariant of the knot up to multiplying by $-1$.
\qed

\subsection{An Analogue of Turaev comultiplication for pseudo-classical curves on $\Sigma$}
\label{sec:DeltaForCurves}

Using the same approach,
we can define a flat version of $\Delta$,
i.e., an analogue of Turaev comultiplication for not knot diagrams
but for oriented pseudo-classical curves without any over/under information in its self-intersections.
More precisely, we will define $\bar{\Delta}$ --- an analogue of Turaev comultiplication for
diagrams of pseudo-classical knots on $\Sigma$
in whose crossings over/under information is ignored.
Unfortunately, in this case, we are forced to use another ground ring.

Let $\bar{\ring}$ be the quotient $\ring \slash \ring_2$,
where $\ring_2$ is the ring whose elements are complex numbers of the form $2a +b i, a,b \in \mathbb{Z}$.
Then we define an algebra $\bar{\algA}$ by analogy with the definition of $\algA$ (see Section~3.3)
with two differences:
\begin{itemize}
\item as a ground ring, we use $\bar{\ring}$ instead of $\ring$,
\item this time, we use the symmetric tensor product of the generators.
\end{itemize}
Let $\bar{\algA}_n$ denote the homogeneous degree $n$ component of $\bar{\algA}$.

Set
$$
\bar{\Delta}(D) =\sum_{x \in \#D}\FlatOneSign(x) [l_1(x)] \otimes [l_2(x)] \in \bar{\algA}_2,
\eqno(4)
$$
where $D \subset \Sigma$ is a pseudo-classical diagram $K$
and $\FlatOneSign(x)$ is the sign of the flat self-intersection point defined in Section~2.3.2.
It is necessary to emphasize that,
in the flat case, we cannot uniquely determine the order of loops corresponding to a crossing,
and this is the reason we use the symmetric tensor product.
The reason we use $\mod{2}$ addition will be clear
from the proof of the next theorem.

\begin{theorem}
\label{theorem:FlatDelta}
If $D_1$ and $D_2$ are pseudo-classical diagrams, 
so that one can be transformed into the other by a finite sequence of flat Reidemeister moves,
then
$\bar{\Delta}(D_1)$ and $\bar{\Delta}(D_2)$
coincide up to simultaneous complex conjugation of all the coefficients. 
\end{theorem}

\emph{Proof}.

The proof coincides with the one of Theorem~\ref{theorem:Delta}
except for two differences.

The first one concerns the invariance under the move $R_2$.
Now, we deal with $\FlatOneSign$ of crossings
for which (see Remark~\ref{remark:FlatBigonSelfIntersection})
$\FlatOneSign(x_1) =\overline{\FlatOneSign(x_2)}$,
where $x_1$ and $x_2$ are crossings involved in flat $R_2$-move.
Hence, the terms on the right-hand side of~(4)
corresponding to $x_1$ and $x_2$ cancel out
in both of the possible situations:
if $\FlatOneSign(x_1) =\pm i$ then we have $i-i =0$,
if $\FlatOneSign(x_1) =1$ we have $1+1 =0 \mod{2}$.

The second difference concerns the relabeling of the diagram,
which, this time,  leads to simultaneous relabeling of intersecting loops.
In the case of flat self-intersection point (see Proposition~\ref{proposition:FlatSignSelfIntersection}),
it implies the complex conjugation of all the coefficients in the right-hand side of~(4).
\qed

\section{Some secondary invariants}
\label{sec:SecondaryInvariants}

\subsection{Comultiplication for unoriented pseudo-classical knots}
\label{sec:Unoriented}

To avoid the dependence of $\Delta(K)$ on the orientation of the knot
we use the following  standard trick:
$$
\Delta_0(|K|) =\Delta(K) +\Delta(-K) \in \algA_2,
$$
where $|K|$ denotes a pseudo-classical knot $K$ with orientation ignored
and $-K$ denotes $K$ with the orientation reversed and the labeling permuted.

\begin{lemma}
\label{lemma:DeltaForUnoriented}
The value $\Delta_0(|K|)$ is an isotopy invariant of the unoriented pseudo-classical knot $K$
up to multiplying by $-1$.
\end{lemma}

\emph{Proof}.

Choose an arbitrary orientation of the knot $K$
and fix a labeling of its diagram $D$.
By Proposition~\ref{proposition:Sign}(2),
the signs of the corresponding crossing in the diagrams $D$ and $-D$ coincide.
By the rule of numbering loops in a crossing,
the first one goes away from the crossing along the overgoing arc in the positive direction,
hence in each crossing of $D$ and $-D$  the loops are permuted
(and their orientations are reversed).
Thus,
$$
\Delta_0(|K|) =\sum_{x \in \#D}\sign(x) \big( [l_1(x)] \otimes [l_2(x)] +[-l_2(x)] \otimes [-l_1(x)] \big).
$$
The sum on the right-hand side coincides (up to the permuting terms in parentheses)
with the one corresponding  to the knot $-K$.

The value $\Delta_0$ is an invariant up to multiplying by $-1$
because the relabeling of the diagram $D$ 
implies (see Proposition~\ref{proposition:Sign}(1)\,)
multiplying the resulting value by $-1$.
\qed

\subsection{A homological consequence of $\Delta$}
\label{sec:HomologicalComultiplication}

Sometimes, it is more convenient to work with homology classes of loops rather than their homotopy classes.
The following invariant of pseudo-classical knots is a homology version of our analogue of Turaev comultiplication:
$$
\Delta_H(K) =\sum_{x \in \#D_{0,H}}\sign(x) [l_1(x)]_H \otimes [l_2(x)]_H,
$$
where $D$ is a diagram of oriented pseudo-classical knot $K$,
$[\cdot]_H \in H_1(\Sigma)$ denotes the homology class of the specified loop,
$x \in \#D_{0,H} \subset \#D$ if and only if
both of $l_1(x),l_2(x)$ are not null-homologous.
In the case of a thickened orientable surface,
an invariant analogous to $\Delta_H$ studied in~\cite{Tarkaev2020_1}.

\subsection{The composition of $\Delta$ and the analogue of Goldman bracket}
\label{sec:CompositionOfDeltaAndGoldman}

Let $K$ be a pseudo-classical knot
represented by a diagram $D$.
Consider the composition of $\Delta$ and the analogue of Goldman bracket (see Section~3.3):
$$
[\Delta](K) =\sum_{x \in \#D}\sign(x) [ l_1(x), l_2(x) ] \in \algA_1.
$$
Note that:

1. Each loop of the form
$[l_1(x), l_2(x) ]$ has at least two crossing less than $D$,
i.e., we obtain the sum of homotopy classes of loops, each of which is simpler than the initial one.

2. As in the previous section, instead of the homotopy classes, we can consider homology classes of the loops
$[l_1(x), l_2(x)]$.
The resulting invariant takes value in the free $\mathbb{Z}$-module generated by elements of the group $H_1(\Sigma)$.

\subsection{The composition of $\bar{\Delta}$ with itself}
\label{sec:barbarDelta}

Let $C \subset \Sigma$ be a generic oriented pseudo-classical curve
given by forgetting the over/under information in a pseudo-classical diagram $D \subset \Sigma$.

Consider
\begin{equation*}
\label{eq:barbarDelta}
\barbar{\Delta}(C) =\sum_{n \geq 1} \barbar{\Delta}_n(C) \in \bar{\algA},
\end{equation*}
where $\barbar{\Delta}_n(C) \in \bar{\algA}_n$
and there is $n_0 \geq 1$ so that $\barbar{\Delta}_n(C) =0$ if $n > n_0$.

The procedure giving the value of $\barbar{\Delta}(C)$ can be described as follows.
$$
\barbar{\Delta}_1(C) =[C], \quad \barbar{\Delta}_2(C) =\bar{\Delta}(C).
$$
To obtain $\barbar{\Delta}_3(C)$ we apply $\bar{\Delta}$ sequentially to each factor in each term of $\barbar{\Delta}_2(C)$:
$$
\barbar{\Delta}_3 =\sum_{x \in \#D} \FlatOneSign(x) \big(\bar{\Delta}(l_1(x)) \otimes [l_2(x)] +[l_1(x)] \otimes \bar{\Delta}(l_2(x)) \big).
$$
Here we set that $\bar{\Delta}(l) =0$ if $l$ is not a pseudo-classical loop.
The resulting sum consists of terms,
each of which is a product of three factors
with a coefficient of the form $\FlatOneSign(x)~\cdot~\FlatOneSign(y)$,
where ``$\cdot$'' stands for the multiplication in the ring $\bar{\ring}$.
Then we analogously obtain $\barbar{\Delta}_4(C)$ from $\barbar{\Delta}_3(C)$,
namely, we apply $\bar{\Delta}$ sequentially to each factors in each term.
As a result, we obtain the sum of products of four homotopy classes of loops
with coefficients of the form $\FlatOneSign(x) \cdot \FlatOneSign(y) \cdot \FlatOneSign(z)$.

And so on: applying $\bar{\Delta}$ we obtain $\barbar{\Delta}_{n+1}(C)$ from $\barbar{\Delta}_n(C)$.

The above process is necessarily finite because applying $\bar{\Delta}$ decreases the number of crossing by at least one.

Note that loops involved in the resulting sum
either are not pseudo-classical
or are homotopic to simple loops
because, otherwise, further applying of $\bar{\Delta}$ is possible.

The above construction
works for knots in a thickened orientable surface as well.
Moreover, in the orientable case, it becomes simpler and stronger,
because we do not need the $\mod{2}$ addition,
since the crossings involved in the flat $R_2$-move are always $+1$ and $-1$.
Therefore, we can use $\mathbb{Z}$ as a ground ring of the corresponding algebra.

\subsection{An analogue of the affine index polynomial}
\label{sec:AnalogueOfAffineIndexPolynomial}

Consider a linear mapping
$f: \algA_2 \to \ring[x^{\pm 1},y^{\pm 1}]$
that sends $[l_1] \otimes [l_2]$ to $x^{\myRe(l_1 \cdot l_2)} y^{\myIm(l_1 \cdot l_2)}$,
 where $l_1,l_2$ are loops in a pseudo-classical diagram $D$
and $l_1 \cdot l_2$ is understood in the sense of Remark~\ref{rem:IntersectionNumberForSubloops}.
Composing $\Delta$ and $f$ we obtain
a polynomial invariant of pseudo-classical knot
$$
P(K) =\sum_{x \in \#D}\sign(x) x^{\myRe(l_1(x) \cdot l_2(x))} y^{\myIm( l_1(x) \cdot l_2(x))},
\eqno(5)
$$
where $D$ is a diagram of the pseudo-classical knot $K$.

Note that the relabeling of $D$ leads to the simultaneous multiplying by $-1$
both signs of all crossings and all intersection numbers involved,
hence the polynomial $P(K)$ is an isotopy invariant of $K$ up to such a transformation.
In the classical situation, when we consider knots in a thickened orientable surface,
the transformation of the affine index polynomial corresponds to changing of the orientation of the surface.

The above definition is slightly different from the definition of the affine index polynomial~(1).
Namely, in~(1) 
all terms have the form $\sign(x) (t^n -1)$.
The reason for using ``$-1$'' is to guarantee the invariance under the first Reidemeister move:
the term corresponding to the crossing involved in the move equal $\sign(x) (t^0 -1) =0$.
But as a result, any crossing $x$ so that $l_1(x) \cdot l_2(x) =0$ does not contribute the resulting value
even when both involved loops are not contractible.
While in~(5)
a crossing does not matter only if at least one of $l_1(x),l_2(x)$ is contractible.



\end{document}